\documentclass[12pt,a4paper,eqno]{amsart}
\usepackage{graphicx,epsfig}
\usepackage{amssymb}

\def\ind{\pmb{1}}

\def\qed{\hfill{$ \Box $}}

\newenvironment{namelist}[1]{%
\begin{list}{}
     {
      
      \settowidth{\labelwidth}{#1}
      \setlength{\leftmargin}{1.1\labelwidth}
               }
      }{%
\end{list}}

\newtheorem{thm}{Theorem}[section]
\newtheorem{cor}[thm]{Corollary}
\newtheorem{lem}[thm]{Lemma}
\newtheorem{prop}[thm]{Proposition}
\theoremstyle{definition}
\newtheorem{defn}[thm]{Definition}
\theoremstyle{remark}
\newtheorem{rem}[thm]{Remark}
\newtheorem{ex}{\large\bf Example}

\begin{document}

\title[Kendall random walks]%
      {Kendall random walks}
\author{Barbara H. Jasiulis-Go{\l}dyn}
           \thanks{ Institute of Mathematics University of
           Wroc{\l}aw, pl. Grunwaldzki 2/4, 50-384 Wroc{\l}aw,
           Poland, e-mail: jasiulis@math.uni.wroc.pl 
           }
\maketitle
%\begin{abstract}

\begin{abstract} The paper deals with a new class of random walks strictly connected with the Pareto distribution.
We consider stochastic processes in the sense of generalized convolution or weak generalized convolution following the idea given in \cite{WeakLevyProc}. The processes are Markov processes in the usual sense. Their structure is similar to perpetuity or autoregressive model. We prove theorem, which describes the magnitude of the fluctuations of random walks generated by generalized convolutions.

We give a construction and basic properties of random walks with respect to the Kendall convolution. We show that they are not classical L\'evy processes. The paper proposes a new technique to cumulate the Pareto-type distributions using a modification of the Williamson transform and contains many new properties of weakly stable probability measure connected with the Kendall convolution. It seems that the Kendall convolution produces a new classes of heavy tailed distributions of Pareto-type.
\end{abstract}

\section{Introduction}

In 2009 Nguyen Van Thu in \cite{Thu}, considering only the Kingman convolution, showed that each generalized convolution, together with an infinitely divisible distribution with respect to this convolution, defined a Markov process, which could be treated as a L\'evy process in the sense of this convolution. The most important example is the Bessel process defined by the Kingman convolution (sometimes called also Bessel convolution) - widely applied and intensively studied in many different areas of mathematics.

The L\'evy processes with respect to generalized convolutions and weak generalized convolutions were introduced in \cite{WeakLevyProc}. The paper deals with the L\'evy processes constructed as the Markov processes in the usual sense. 

We consider here Markov chain $\{X_n \colon n \in \mathbb{N}_0 \}$, where $\mathbb{N}_0 = \mathbb{N} \cup \{0\}$, given by the construction proposed in \cite{WeakLevyProc}. 

The possibility of many different interpretations of the cumulation of independent steps $\Delta X_1 ,\dots, \Delta X_k$ into the state $X_k$-cumulation, which is not just simple addition, promises wide applicability of such processes.
Similarly to the classical theory, the discrete time processes can be built on the basis of any step distribution (in the sense that distribution of unit step does not have to be infinitely divisible in any sense). The main generalization is that instead of classical addition we consider here two binary operations. The first one is the generalized convolution, which is an associative and commutative operation on the set of probability measures on the Borel subsets of the positive half line. We also consider generalized convolution extended to the Borel subsets of the real line in the sense given in \cite{jas}, in particular weak generalized convolution. Since weakly stable probability measure, which is strictly connected with the weak generalized convolution, is a natural generalization of symmetric $\alpha-$ stable distribution, it is worth investigating and it can be used in applications.\\

The next section deals with the definition and the existence of the Markov chains based on generalized convolutions and weak generalized convolutions. We recall construction given in \cite{WeakLevyProc}. A theorem describing the magnitude of the fluctuations of constructed random walks will be proved.

In the third section we study properties of random walks under the Kendall generalized and weak Kendall generalized convolutions. If we consider the Kendall convolution for the probability measures on the positive half line, then the main tool which we use is a transform called homomorphism:
$$
h(\delta_t)=(1-t^{\alpha})_+
$$
for $\alpha > 0$. For weakly stable probability measure with the characteristic function:
$$
\widehat{\mu}_{\alpha}(t)=(1-|t|^{\alpha})_+,
$$
for $0 < \alpha \leqslant 1$, we arrive at random walks on the real line.

In the Kendall generalized convolution case we obtain the following random walk:
$$
X_{n+1} = \left(X_n \vee \Delta X_{n+1} \right) \theta_n^{Q_n} \quad a.e.,
$$
where $\vee$ denotes maximum, $(\theta_n)$ is proper i.i.d. sequence with the Pareto distribution with density $\pi_{2\alpha}(dx)=\frac{2\alpha}{x^{2\alpha+1}} \mathbf{1}_{(1,\infty)} dx$ such that $\theta_n$ is independent of $\left(X_n \vee \Delta X_{n+1} \right)$. The random variables $(Q_n)$ take values $0$ and $1$ provided that we know the position of $\left(X_n \vee \Delta X_{n+1} \right)$.

We show that the obtained random walks are not the L\'evy processes in the usual sense.
Their structure is similar to the first order random coefficients autoregressive model (see e.g. \cite{Darek}) but with a different dependence structure.
In order to get properties of constructed stochastic processes we consider here two types of random walks: the first one under generalized convolution and weak generalized convolution and, the second one: random walks associated with them, which are stochastic processes in the usual sense (under the classical convolution). Usually some assumptions on the associated random walks imply properties of random walks in the generalized convolution sense. For instance, we prove an analogue of the law of iterated logarithm for constructed random walk under assumptions about associated random walk. \\

We also present many new properties of the weakly stable probability measure $\mu_{\alpha}$, in particular we obtain a characterization of the Pareto distribution $\pi_{2\alpha}$ in terms of $\mu_\alpha$ for $0 < \alpha \leqslant 1$. We also obtain the connection between $\mu_{\alpha}$ and $\mu_1$, which is really useful because instead of investigating random walk under the Kendall convolution generated by $\mu_{\alpha}$, it is sufficient to consider the case $\alpha =1$, which is easier. These properties are also useful in the study of associated random walks, i.e. classical random walks with unit step distribution, which has the factor $\mu_{\alpha}$.
\\

In order to start our constructions we need some basic facts on generalized convolution and weak generalized convolution:

The main mathematical tool used here, is the generalized convolution defined by Urbanik (see \cite{Urbanik64}) on the set $\mathcal{P}_+$ of probability measures on the Borel subsets of the positive half line. For simplicity we will use notation $T_a$ for the rescaling operator defined by $(T_a\lambda)(A) = \lambda ({A/a})$ for every Borel set $A$ when $a \neq 0$, and $T_0 \lambda = \delta_0$ is the probability measure concentrated at zero.

\begin{defn}\label{df:1}
A commutative and associative $\mathcal{P}_{+}$-valued binary
operation $\diamond$ defined on $\mathcal{P}_{+}^2$ is called a
\emph{generalized convolution} if for all $\lambda,\lambda_1,\lambda_2 \in \mathcal{P}_{+}$ and $a \geqslant 0$ we have:
\begin{namelist}{llll}
\item[{\rm (i)}] $\delta_0 \diamond \lambda = \lambda$ ; %
\item[{\rm (ii)}] $(p\lambda_1 +(1-p)\lambda_2) \diamond \lambda = p ( \lambda_1 \diamond \lambda) + (1-p)(\lambda_2
          \diamond \lambda) $ whenever $p\in [0,1]$; %
\item[{\rm (iii)}] $T_a(\lambda_1 \diamond \lambda_2) = (T_a \lambda_1) \diamond (T_a \lambda_2)$ ; %
\item[{\rm (iv)}] if $\lambda_n \rightarrow \lambda$ then $\lambda_n \diamond \eta \rightarrow \lambda \diamond
           \eta$  for all $\eta \in \mathcal{P}_{+}$ and $\lambda_n \in \mathcal{P}_{+}$; %
\item[{\rm (v)}] there exists a sequence $(c_n)_{n\in\mathbb{N}}$ of positive numbers such that the sequence $T_{c_n}
           \delta_1^{\diamond n}$ converges to a probability measure different from $\delta_0$,
\end{namelist}
where $\rightarrow$ denotes a weak convergence of probability measures.
\end{defn}

We call the set $(\mathcal{P}_+, \diamond)$ \emph{generalized convolution algebra}. A continuous mapping $h: \mathcal{P}_+ \to \mathbb{R}$ such that
\begin{itemize}
\item $h(p\lambda + (1-p)\nu)= p h(\lambda) + (1-p) h(\nu)$,
\item $h(\lambda  \diamond \nu)=h(\lambda) h(\nu)$
\end{itemize}
for all $\lambda, \nu \in \mathcal{P}_+$ and $p \in (0,1)$ is called \emph{homomorphism of $(\mathcal{P}_+, \diamond)$}. For every probability measure $\lambda \in \mathcal{P}_+$ we have
$$
h(\lambda) = \int_{\mathbb{R}_+} h(\delta_x) \lambda(dx).
$$
The algebra $(\mathcal{P}_+, \diamond)$ is regular if it admits a nontrivial homomorphism, i.e. such h that $h \not\equiv 0$ and $h \not\equiv 1$.
Since for all $\lambda_1,\lambda_2\in\mathcal{P}_+$ we have
    $$
    \lambda_1 \diamond \lambda_2(A)=\int_{\mathbb{R}_+} \int_{\mathbb{R}_+} \rho_{x,y}(A) \lambda_1(dx) \lambda_2(dy),
    $$
every generalized convolution is uniquely determined by the {\em probability kernel}
$$
\rho_{x,y} := \delta_x \diamond \delta_y.
$$
Evidently
\begin{itemize}
    \item
    $ \rho_{x,0}=\delta_x$,
    \item
     $ \rho_{x,y}=\rho_{y,x}$,
    \item
    $ T_c\rho_{x,y}=\rho_{cx,cy}$,
    \item
    $ \rho_{x,y}=T_v\rho_{z,1}$, where $v=x\vee y,\; z=\frac{x\wedge y}{x\vee y}$
        \end{itemize}
for all $x,y,c \geqslant 0$. The origin of the generalized convolution can be found in the Kingman paper (\cite{King}), where the first example of random walk under generalized convolution (called the Kingman or Bessel convolution) is considered. The Kingman convolution has a natural interpretation at the interference phenomena (see \cite{plastimo}). In series of papers \cite{KU1a}--\cite{Urbanik64} Urbanik developed the theory of generalized convolutions. In \cite{vol1}--\cite{vol5} Vol'kovich was investigating this theory. Many open problems connected with generalized convolutions were given by Vol'kovich, Toledano-Ketai and Avros (see \cite{vol2}), Hazod (\cite{Hazod}) or Van Thu (\cite{Thu}). In 2012 there appeared \cite{JasKula} about generalized convolutions in the non-commutative probability theory. Hence, it can be supposed that the results obtained here will be also useful in the non-commutative probability theory. On the other hand, the general theory of the L\'evy processes in the generalized convolution sense was first established in \cite{WeakLevyProc}. In this paper we develop that theory giving explicite recipe for random walks with respect to the Kendall convolution, presenting new results on this topic and proposing a new technique of the Williamson transform to investigate constructed stochastic processes.

\begin{ex}The best known examples of generalized convolutions are the following:
\begin{namelist}{ll}
\item[$\bullet$]
$\alpha-$convolution $\star_{\alpha}$ for $\alpha>0$ given by the formula
$\delta_a\star_{\alpha}\delta_b=\delta_c$, where
$c^{\alpha}=a^{\alpha}+b^{\alpha} $; the corresponding homomorphism
$h(\delta_t)=\exp\{-t^{\alpha}\}$,
\item[$\bullet$]
max-convolution $\star_{\infty}$ given by
$\delta_a\star_{\infty}\delta_b=\delta_c$, where $c=a\vee b$; the corresponding homomorphism $h(\delta_t)=\ind_{[0,1]}(t)$,
\item[$\bullet$]
symmetric convolution $\star_{1,1}$ given by
$\delta_a\star_{1,1}\delta_b=\frac{1}{2}\delta_{a+b}+
\frac{1}{2}\delta_{|a-b|}$ ; the corresponding homomorphism $h(\delta_t)=cos(t)$,
\item[$\bullet$] the Kendall convolution $\triangle_{\alpha}$ for $\alpha>0$ with the probability kernel
$\delta_a\triangle_{\alpha}\delta_b=\left(1-\frac{a^{\alpha}}{b^{\alpha}} \right)\delta_b+\frac{a^{\alpha}}{b^{\alpha}}T_b\pi_{2\alpha}$
for $0\leqslant a \leqslant b$, where $\pi_{2\alpha}$ is the Pareto distribution with the density
$\pi_{2\alpha}(dx)=\frac{2\alpha}{x^{2\alpha+1}}\ind_{(1,\infty)}(x)dx$ and $h(\delta_t)=(1-t^{\alpha})_+$,
\item[$\bullet$]
the Kingman convolution $\star_{1,\beta}$ $(\beta>1):$
 \begin{eqnarray*}
  \lefteqn{
\delta_a\star_{1,\beta}\delta_b(dx) = B\left(\frac{\beta-1}{2}, \frac{1}{2}\right) } \\
&  & \cdot \frac{\left[ \left(x^2-(a-b)^2\right) \left((a+b)^2-x^2\right)\right]^
{\frac{\beta-3}{2}}} {(2ab)^{\beta-3}}\ind_{[|a-b|,a+b]}(x)dx
 \end{eqnarray*}
and
$$
h(\delta_t)=\Gamma\left(\frac{\beta}{2}\right) \left(\frac{2}{t}\right)^{\frac{\beta}{2}-1} J_{\frac{\beta}{2}-1}(t) ,
$$
where $J_{\frac{\beta}{2}-1}(t)$ is the Bessel function of order $\frac{\beta}{2}-1$ and $B(a,b)$ is the Beta function with parameters $a$ and $b$,
\item[$\bullet$]
the Kucharczak convolution $\circ_{\alpha}$ $(0<\alpha<1)$:
$$
\delta_a\circ_{\alpha}\delta_b(dx)
=\frac{ sin(\pi\alpha) \cdot (ab)^{\alpha} (2x-a-b)}{\pi \cdot (x(x-a-b))^{\alpha}\cdot (x-a) \cdot (x-b)}\pmb{1}_{[(a^{\alpha}+b^{\alpha})^{1/\alpha}, \infty)}(x)dx
$$
and $h(\delta_t)(t)=\Gamma(\alpha)^{-1}\Gamma(\alpha,t)$, where $\Gamma(\alpha,t)$ is the incomplete $\Gamma-$ function; for more details see \cite{KuU2}.
\end{namelist}
\end{ex}\label{ex:1}
In \cite{jas} one can find the definition and basic properties of generalized convolution on set $\mathcal{P}$ of probability measures on the Borel subsets of the real line. An example of such object is a weak generalized convolution connected with weakly stable distribution. We recall the definition of weak stability of probability measures.
\begin{defn}\label{df:2}
A probability measure $\mu \in \mathcal{P}$ is \emph{weakly stable} if
\[
\forall \; a,b \in \mathbb{R} \; \exists \lambda \in \mathcal{P} \quad  T_a \mu * T_b \mu = \mu \circ \lambda
\]
or equivalently
\[
\forall \; \lambda_1,\lambda_2 \in \mathcal{P}\; \exists \lambda \in \mathcal{P} \quad  \mu \circ \lambda_1 * \mu \circ \lambda_2  = \mu \circ \lambda,
\]
\end{defn}

where $(\mu \circ \lambda)(A) = \int \mu({A/s}) \lambda(ds)$ for every Borel set $A$. From Theorem 6 in \cite{MOU} we know that if $\mu$ is a weakly stable probability measure on a separable Banach space $\mathbb{E}$, then either there exists $a\in \mathbb{E}$ such that $\mu = \delta_a$ or there exists $a \in \mathbb{E} \setminus \{0\}$ such that $\mu = \frac{1}{2} (\delta_a + \delta_{-a})$ or $\mu(\{a\}) = 0$ for every $a \in \mathbb{E}$. In this paper we consider a non-trivial weakly stable measures, i.e. weakly stable probability measures without any atoms.
In \cite{Mis2} Misiewicz defined a weak generalized convolution:
\begin{defn}\label{df:3}
For weakly stable measure $\mu$ a \emph{weak generalized convolution} of $\lambda_1$ and
$\lambda_2$ (notation $\lambda_1 \otimes_{\mu} \lambda_2$) is defined by the following formula:
$$
  \lambda_1 \otimes_{\mu} \lambda_2=
  \left\{\begin{array}{ll}
  \lambda & if \quad \mu \quad is \quad non-symmetric, \\
  \left|\lambda\right| & if \quad\mu\quad is \quad symmetric,
            \end{array}\right.
$$
\end{defn}
where $|\lambda|= \mathcal{L}(|\theta|)$ if $\lambda= \mathcal{L}(\theta)$. Instead of $|\lambda|\in\mathcal{P}_+$ we can take its symmetrization $\widetilde{\lambda} =\frac{1}{2} |\lambda|+\frac{1}{2}T_{-1}|\lambda|$, since in both cases we have uniqueness of $\lambda_1 \otimes_{\mu} \lambda_2$. We will consider only the case of weak generalized convolution defined by  $\widetilde{\lambda}$, which is more convenient.\\
It can be shown that for each weakly stable probability measure $\mu$ weak generalized convolution $\otimes_{\mu}$ has properties (i)-(iv) on $\mathcal{P}$ but it does not have to have condition (v).
A wide discussion on condition (v) for generalized convolution on $\mathcal{P}$ can be found in \cite{jas}. In particular if a weakly stable measure belongs to the domain of attraction of some strictly stable measure, then it generates weak generalized convolution having property (v).
Similarly to the generalized convolution theory every weak generalized convolution $\otimes_{\mu}$ is uniquely determined by a \emph{ weak probability kernel}:
$$ 
\rho_{z,1}=\delta_z \otimes_{\mu} \delta_1,
$$
since $ \rho_{x,y}=T_v\rho_{z,1}$ for $v=|x|\vee|y|,\; z=\frac{|x|\wedge |y|}{|x|\vee |y|}$ for all $x,y\in\mathbb{R}$, where $\vee$ and $\wedge$ denote maximum and minimum respectively. Moreover,
$$\lambda_1 \otimes_{\mu} \lambda_2(A)=\int_{\mathbb{R}} \int_{\mathbb{R}} \rho_{x,y}(A) \lambda_1(dx) \lambda_2(dy)$$
for all $\lambda_1,\lambda_2 \in \mathcal{P}$. It follows that it is sufficient to define $\rho_{z,1}$ for $|z|\leqslant 1$.

\vspace{2mm}

\begin{ex}
The most popular examples of weakly stable distributions are the following:
\begin{namelist}{ll}
\item[$\bullet$] symmetric $\alpha-$ stable measure $\gamma_{\alpha}$ with the characteristic function $\widehat\gamma_{\alpha}(t)= \exp\{- A|t|^{\alpha} \}_+$, where $A$ is a constant, for $0<\alpha\leqslant 2$, defines weak generalized convolution $\otimes_{\gamma_{\alpha}}$ by
    $$
    \delta_a \otimes_{\gamma_{\alpha}} \delta_b=\widetilde{\delta}_c,
    $$
    where
    $|c|^{\alpha}=|a|^{\alpha}+|b|^{\alpha}$, since
		$$
		T_a \gamma_{\alpha} * T_b \gamma_{\alpha} = T_c \gamma_{\alpha};
		$$
\item[$\bullet$]  probability measure with the characteristic function $\widehat{\mu}_{\alpha}(t)=( 1 - |t|^{\alpha} )_+$ generates the Kendall convolution $\otimes_{\mu_{\alpha}} (0<\alpha\leqslant 1):$
    $$
    \delta_a\otimes_{\mu_{\alpha}} \delta_b=\left(1-\frac{|a|^{\alpha}}{|b|^{\alpha}}\right)\widetilde{\delta}_b+ \frac{|a|^{\alpha}}{|b|^{\alpha}}T_b\widetilde{\pi}_{2\alpha},
    $$
    for $|a|<|b|$, where $\widetilde{\pi}_{2\alpha}$ is the symmetrization of the Pareto distribution with the density
    $\widetilde{\pi}_{2\alpha}(dx)=\frac{\alpha}{|x|^{2\alpha+1}} \ind_{(1,\infty)}(|x|)dx$,
		since
		$$
		T_a \mu_{\alpha} * T_b \mu_{\alpha} = T_c \mu_{\alpha}
		$$
		for all $a,b \in \mathbb{R}$;
\item[$\bullet$]  uniform distribution on the unit sphere $S_{n-1}\subset \mathbb{R}^n$ corresponding to the Kingman convolution,
\item[$\bullet$]  every k-dimensional projection of a weakly stable random vector is weakly stable,
\item[$\bullet$]  distributions introduced by Cambanis, Keener and Simons in \cite{cks}.
\end{namelist}
\end{ex}
Notice that we can produce new classes of weak generalized convolutions generated by generalized convolution. For details see \cite{JasKula} or \cite{KuU2}.
The weak generalized convolution theory contains many open problems. Jarczyk and Misiewicz (\cite{JarMis}) consider pseudoisotropic distributions connected with weak stability and solve many open problems from the point of view of functional equations. In \cite{Grazka} Mazurkiewicz describes weakly stable distributions connected with distributions introduced by Cambanis, Keener and Simons. The recipe for weak generalized convolution connected with these distributions is still an open problem.

\section{The Markov chain under generalized and weak generalized convolution}

Following paper \cite{WeakLevyProc} we consider here a discrete time L\'evy process (independent increments random walk) in the sense of generalized convolution.
More precisely: we will investigate $\{X_n: n \in \mathbb{N}_0\}$, where $X_0 = 0$ and $\mathbb{N}_0=\mathbb{N}\cup\{0\}$, based on the set of i.i.d. random variables $\{\Delta X_k: k\in\mathbb{N} \}$ with distribution $\nu \in \mathcal{P}_+$. Random element $X_n$ is a kind of cummulation of variables $\Delta X_1, \Delta X_2, \cdots, \Delta X_n$ in the following sense:
$$
\mathcal{L}(X_n) = \mathcal{L}(\Delta X_{1}) \diamond \cdots \diamond \mathcal{L}(\Delta X_{n})
$$
for all $n\geqslant 1$. Moreover, we assume that the increment of the process from $X_n$ to $X_{n+k}$ in the sense of generalized convolution depends only on $\Delta X_{n+1}, \cdots, \Delta X_{n+k}$, i.e.
$$
\mathcal{L}(X_n) \diamond \mathcal{L}(\Delta X_{n+1}) \diamond \cdots \diamond \mathcal{L}(\Delta X_{n+k}) = \mathcal{L}(X_{n+k}).
$$
for all $n,k \geqslant 1$.
The existence theorem of the process $\{X_n: n \in \mathbb{N}_0\}$ was proved in \cite{WeakLevyProc}.

\begin{thm}\label{thm:1}
There exists a Markov process $\{X_n: n \in \mathbb{N}_0\}$ with the transition probabilities
$$
P_{k,n}(x,A)=P(X_n\in A|X_k=x):= \delta_x \diamond \nu^{\diamond n-k}(A),
$$
where $x\geqslant 0$, $n,k \in \mathbb{N}$ and $A\in\mathcal{B}((0,\infty))$.
\end{thm}

\vspace{5ex}
On the other hand we can consider stochastic processes generated by a weakly stable probability measure 
$\mu$ (in the sense of weak generalized convolution $\otimes_{\mu}$). Then we construct two associated processes $\{\widetilde{X}_n: n\in \mathbb{N}_0\}$ and $\{\widetilde{S}_n: n\in \mathbb{N}_0\}$ such that $\{\widetilde{X}_n: n\in \mathbb{N}_0\}$ is defined as above process $\{X_n: n\in \mathbb{N}_0\}$  but on the real line. We assume, without loss of generality, that $\widetilde{X}_0=0$. This means that the families $\{\Delta \widetilde{X}_k\}$ and $\{\Delta \widetilde{X}_{k,n}\}$ are also specified such that $(\Delta \widetilde{X}_k)_{k\in\mathbb{N}}$ are independent identically distributed with distribution $\nu\in \mathcal{P}$ and $\widetilde{X}_k \perp \Delta \widetilde{X}_{k,n}$ for every $k\in\mathbb{N}$ and $n>k$. Additionally, for a sequence  $Y, (Y_i)_{i\in\mathbb{N}_0}$ of i.i.d. random variables with weakly stable distribution $\mu$ we define
$$
\widetilde{S}_n:=\sum_{k=1}^n  \left( \Delta \widetilde{X}_k \cdot Y_k \right),
$$
where $(Y_i)_{i\in\mathbb{N}_0}$ and $(\Delta \widetilde{X}_i)_{i\in\mathbb{N}}$ are independent.
Since $\widetilde{X}_0=0$, we have $S_0=0$. Notice that $\{\widetilde{S}_n: n\in \mathbb{N}_0\}$ is a discrete time L\'evy process in the classical sense such that
$\widetilde{S}_n \stackrel{d}{=} \widetilde{X}_n Y$.
\\
In much the same way as in the previous one, for all $k,n\in\mathbb{N}, \; k\leqslant n$ and the Borel set $A\in \mathcal{B}(\mathbb{R})$, we see that the transition probabilities are given by
$$
\widetilde{P}_{k,n}(x,A)=P(\widetilde{X}_n\in A|\widetilde{X}_k=x):= \delta_x \otimes_{\mu} \nu^{\otimes_{\mu} n-k}(A), \quad x \in \mathbb{R}.
$$
The proof of the next corollary is basically the same as the proof of Theorem \ref{thm:1} .\\
\begin{cor}\label{cor:1}
Let $\mu$ be a non-trivial weakly stable probability measure.
$\{\widetilde{X}_n \colon n \in \mathbb{N}_0 \}$ is a Markov process.
\end{cor}

Using the \emph{ characterizing exponent}
$$
\varkappa:=\varkappa(\mu)=\sup \biggl\{ p \in [0,2] \colon
\int_{\mathbb{R}}|x|^p \mu(dx) < \infty \biggr\}
$$
for the weakly stable distribution $\mu$ we can prove the following theorem describing the magnitude of the fluctuations of random walk under weak generalized convolution. Parameter $\varkappa(\mu)$ plays a similar role as the parameter $\alpha$ for $\alpha$-stable distribution. More information about characterizing exponent one can find in  \cite{misjas2}.

\begin{thm}\label{thm:2}
Let $\mu$ be a non-trivial weakly stable probability measure with $\varkappa(\mu)>0$ and let $\{\widetilde{X}_n: n\in\mathbb{N}_0\}$ be the random walk, under weak generalized convolution $\otimes_{\mu}$, defined above with i.i.d. increments $(\Delta \widetilde{X}_i)_{i\in\mathbb{N}}$ with distribution $\nu$. Let $Y, Y_1, Y_2, \cdots$ be i.i.d. sequence of random variables with distribution $\mu$, which is independent of $(\Delta \widetilde{X}_i)_{i\in\mathbb{N}}$. If there exist sequences of positive numbers $(a_n)_{n\in\mathbb{N}}$ and $(b_n)_{n\in\mathbb{N}}$ such that $a_n \nearrow \varkappa, \; b_n \to \infty$ and
$$
\limsup_{n\to\infty} b_n^{-1} \mathbb{E}(|\widetilde{S}_n|^{a_n})=c \in (0,\infty),
$$
then for every  sequence of positive numbers $(c_n)_{n\in\mathbb{N}}$ such that
$$
\sum\limits_{n=1}^{\infty} c_n^{-1} < \infty,
$$
we have
$$
\mathbf{P} \left( \bigcup_{n=1}^{\infty} \bigcap_{k=n}^{\infty} \left\{ |\widetilde{X}_k|\leqslant \left( \frac{c_k b_k}{d_k} \right)^{1/a_k} \right\} \right) = 1,
$$
where $d_n=E|Y|^{a_n}$.
\end{thm}
\noindent {\bf Proof.} Let $A_n = \left\{ |\widetilde{X}_n|^{a_n} > \frac{c_n b_n}{d_n} \right\}$ and $a_n \nearrow \varkappa, \; b_n \to \infty$. Since $\widetilde{S}_n:=\sum_{k=1}^n  \left( \Delta \widetilde{X}_k \cdot Y_k \right) \stackrel{d}{=} \widetilde{X}_n Y$ for $\widetilde{X}_n, \; Y$ independent, we have
$$
\mathbb{E}(|\widetilde{S}_n|^{a_n})=\mathbb{E}(|\widetilde{X}_n|^{a_n}) \mathbb{E}(|Y|^{a_n}).
$$
%Moreover $E(|Y_1|^{a_n}) < \infty$ for every $n\in \mathbb{N}_0$. It yields
%$$
%\limsup_{n\to\infty} b_n^{-1} E(|\widetilde{X}_n|^{a_n})= \limsup_{n\to\infty} \frac{c}{ E(|Y_1|^{a_n})} < \infty.
%$$
By the Tchebyshev inequality there exists $n_0 \in \mathbb{N}_0$ such that
$$
\mathbf{P}(A_n) \leqslant \frac{d_n \mathbb{E}(|\widetilde{X}_n|^{a_n})}{b_n c_n} = \frac{ \mathbb{E}(|\widetilde{S}_n|^{a_n})}{b_n c_n} \leqslant \frac{c}{c_n }
$$
for every $n\geqslant n_0$. Since $\sum\limits_{n=1}^{\infty} c_n^{-1} < \infty$ we arrive at
$\sum_{n=1}^{\infty} \mathbf{P}(A_n) < \infty$.
Now it is sufficient to apply the Borel-Cantelli lemma. \qed

\section{Random walk under the Kendall generalized convolution}
In 2011 in \cite{misjas1} the authors showed that  the Kendall convolution for $\delta_x$ and $\delta_1$ is the unique generalized convolution which can be written as a convex linear combination of two fixed probability measures such that only coefficients of this combination depend on x. Moreover it can
be shown (see Theorem 1 in \cite{misjas1})
that if $0<\alpha \leqslant 1$, then the Kendall convolution is a weak generalized convolution with respect to
a symmetric weakly stable
measure $\mu_{\alpha}$ with the density function
$$
\mu_{\alpha}(dy)=\frac{\alpha}{\pi y}\int\limits_0^a sin(ty)t^{\alpha-1}dtdy
$$
and the characteristic function $\widehat{\mu}_{\alpha}(t)=(1-|t|^{\alpha})_+$.

In \cite{Neslehova} one can find connections between the Kendall convolution and the Archimedean copulas theory, i.e. a generator of Archimedean copula is the homomorphism of the Kendall convolution. It follows that we can also use the Williamson transform (see \cite{Williamson}) to get our results. 

For a random variable $X \sim \nu$ with cummulative distribution function $F$, classical Williamson transform is defined by:
$$
\mathcal{M}_dF(t) = \int\limits_t^{\infty} \left(1-\frac{t}{x}\right)^{d-1} dF(x) = \left\{
\begin{array}{ll}
  \mathbb{E} \left(1-\frac{t}{X}\right)^{d-1}_+ & if \quad t > 0, \\
   1-F(0) & if \quad t=0,
            \end{array}\right.
$$
where $d \geqslant 2$ is integer.

In this paper we investigate the random walks under the Kendall convolution using technique of homomorphism (resp. characteristic function) for given generalized convolution (resp. weak generalized convolution), which is strictly connected with a modification of the Williamson transform. To see this notice that homomorphism
$$
h(\delta_t)=(1-t)_+
$$
is the Williamson transform for probability measure $\delta_1$ and $d=2$.

We use a modification of the Williamson transform given by:
\begin{eqnarray*}
\lefteqn{ \Phi_{\nu}(t) = h(T_t \nu) = \int\limits_{0}^{\infty} \left(1-(ts)^{\alpha}\right)_+ \nu(ds) }\\
& = & F\left(\frac{1}{t}\right) - t^{\alpha} \int\limits_{0}^{1/t} s^{\alpha} \nu(ds) = \alpha t^{\alpha} \int\limits_{0}^{1/t} s^{\alpha-1} F(s) (ds),
\end{eqnarray*}
which is easy to invert. Since 
$$
\alpha^{-1} t^{\alpha} \Phi_{\nu}\left(\frac{1}{t}\right) = 
\int\limits_{0}^{t} s^{\alpha-1} F(s) (ds),
$$
we have
$$
F(t) = \Phi_{\nu}\left(\frac{1}{t}\right) + \alpha^{-1} t \; \frac{d}{dt} \left[ \Phi_{\nu}\left(\frac{1}{t}\right) \right].
$$
\vspace{5ex}

This section is devoted to the random walk under the Kendall convolution $\triangle_{\alpha}$, $\alpha>0$, with unit step with distribution $\nu$ concentrated on the positive half line. Next, we extend our construction obtaining the random walk under weak Kendall convolution with unit steps having distribution on the real line. In particular we consider a random walk under the Kendall convolution in case $\nu=\delta_1$ or $\nu=\widetilde{\delta}_1$ under weak Kendall convolution.\\
\\
Let $(\Delta X_i)_{i\in\mathbb{N}}$ be a sequence of independent identically distributed random variables with distribution $\nu \in \mathcal{P}_+$ and let $(\theta_i)_{i\in\mathbb{N}_0}$ be a sequence of $i.i.d.$ random variables with the Pareto distribution $\pi_{2\alpha}$ such that $\theta_n$ is independent of $X_n$ and $\Delta X_{n+1}$ for all $n\in\mathbb{N}$. The Markov chain $\{X_n : n\in\mathbb{N}_0\}$ with $\Delta X_1\sim\nu$ is such that
$$
\lambda_{0,n,\alpha}(\nu):= \mathcal{L}(X_n)=\nu^{\triangle_{\alpha}n}
$$
for $n\in\mathbb{N}_0$ and $\nu^{\triangle_{\alpha}0}=\delta_0$.

Our construction of random walk implies that $X_n$ and $\Delta X_{n,n+k}$ are independent and $X_k\stackrel{d}{=}\Delta X_{n,n+k}$ for all $n,k\geqslant 0$. \\
Without loss of generality we can assume that we start from zero, i.e. $X_0=0$ a.e.  It remains to find the measures $\lambda_{0,n,\alpha}(\nu)$ for $n \geqslant 2$.

\begin{prop}\label{prop:1}
Let $\nu \in \mathcal{P}_+$. For each natural number $n\geqslant 2$ and $\alpha>0$ we have
$$
\lambda_{0,n,\alpha}(\nu)(0,x)= \Phi_{\nu}^n\left( \frac{1}{x} \right) + \alpha^{-1} x \frac{d}{dx} \left(\Phi_{\nu}^n\left( \frac{1}{x} \right)\right)
$$
or equivalently
$$
\lambda_{0,n,\alpha}(\nu)(0,x)= \frac{d}{ds} \left( s \Phi_{\nu}^n\left( s^{-1/\alpha} \right)\right) |_{s=x^{\alpha}},
$$
where
$$
\Phi_{\nu}\left( t \right) = h\left( T_t \nu \right) = \int_0^{\infty} \left(1- (xt)^{\alpha}\right)_+ \nu(dx)
$$
is the homomorphism of the unit step variable $\triangle X_1$.
\end{prop}
\noindent {\bf Proof.}
Since $\mu_{\alpha}$ is weakly stable, then
$$
\left(\mu_{\alpha} \circ \nu \right)^{\ast n} = \mu_{\alpha} \circ \nu^{\otimes_{\mu_{\alpha}} n}.
$$
It follows
$$
\Phi_{\nu}^n\left( t \right) = \int_0^{\infty} \left(1- (xt)^{\alpha}\right)_+ \lambda_{0,n,\alpha}(\nu)(dx),
$$
i.e.
$$
\Phi_{\nu}^n\left( t \right) = \Phi_{\nu^{\triangle_{\alpha} n}}\left( t \right).
$$
Consequently, if we convert the transform, we arrive at
$$
\lambda_{0,n,\alpha}(\nu)(0,x)= \Phi_{\nu}^n\left( \frac{1}{x} \right) + \alpha^{-1} x \frac{d}{dx} \left(\Phi_{\nu}^n\left( \frac{1}{x} \right)\right).
$$
In order to obtain equivalent equation it is sufficient to calculate the derivative.
\qed

\begin{cor}\label{cor:2}
Let $\beta_{a,b}$ be the Beta distribution with parameters $a,b>0$, i.e. the distribution with the density function
$$
\beta_{a,b}(dx)=\frac{\Gamma(a+b)}{\Gamma(a) \Gamma(b)}x^{a-1}(1-x)^{b-1}\pmb{1}_{(0,1)}(x)dx.
$$
Then for each natural number $n\geqslant 2$ and $\alpha>0$ the distribution function of the measure $\lambda_{0,n,\alpha}(\beta_{a,b})$ is given by
$$
\frac{d}{ds} \left[\,s\Bigl(B(a,b,s^{1/\alpha})- \frac{\Gamma(a+ \alpha)\Gamma(a+b)}{s \Gamma(a) \Gamma(a+b+\alpha)}\, B(a+\alpha,b,s^{1/\alpha})\Bigr)^{\!n} \,\right] \bigg|_{s=x^{\alpha}}
$$
for $B(a,b,s)=\int_0^s \frac{\Gamma(a+b)}{\Gamma(a) \Gamma(b)}\,x^{a-1}(1-x)^{b-1}dx$. In particular for the uniform distribution $\mathcal{U}(0,1)$ we have
\begin{eqnarray*}
\lefteqn{
\left(\lambda_{0,n,\alpha}(\beta_{1,1})\right)([0,x))= \left(\frac{\alpha}{\alpha+1}\right)^n\left(1+\frac{n}{\alpha}\right)x^n \pmb{1}_{[0,1)}(x)}\\
& & +\left(1-\frac{1}{(\alpha+1)x^{\alpha}}\right)^{n-1}\left(1+\frac{n-1}{( \alpha+1)x^{\alpha}}\right)\pmb{1}_{[1,\infty)}(x).
\end{eqnarray*}
\end{cor}

\vspace{2mm}

\begin{cor}\label{cor:3}
Let $\gamma_{a,b}$ be the Gamma distribution with parameters $a,b>0$ with the density
$$
\gamma_{a,b}(dx)=\frac{b^a}{\Gamma(a) }x^{a-1}e^{-bx}\pmb{1}_{(0,\infty)}(x)dx.
$$
 Then for each natural number $n\geqslant 2$ and $\alpha>0$ the measure $\lambda_{0,n,\alpha}(\gamma_{a,b})$ has the cumulative distribution function
$$
\frac{d}{ds} \left[\,s\Bigl(\Gamma(a,b,s^{1/\alpha})-\frac{ \Gamma(a+\alpha)}{s\Gamma(a) b^{\alpha}}\Gamma(a+\alpha,b, s^{1/\alpha})\Bigr)^{\!n}\,\right] \bigg|_{s=x^{\alpha}}
$$
for $\Gamma(a,b,s)=\int_0^s \frac{b^a}{\Gamma(a) }x^{a-1}e^{-bx} dx$.
\end{cor}
In case $\nu=\delta_1$ it is clear that $P(X_1=1) =1$ and $X_2\stackrel{d}{=}\theta_1$.
\begin{cor}\label{cor:4}
For each natural number $n\geqslant 2$ and real number $\alpha>0$
$$
\left(\lambda_{0,n,\alpha}(\delta_1)\right)(dx)=\frac{\alpha n (n-1)}{x^{2\alpha+1}}\left(1-\frac{1}{x^{\alpha}}\right)^{n-2} \pmb{1}_{[1,\infty)}(x)dx.
$$
\end{cor}
Notice that
$$
\lambda_{0,2,\alpha}(\delta_1) = \pi_{2\alpha}
$$
is the Pareto distribution and by weak stability of $\mu_{\alpha}$ we have
$$
\mu_{\alpha} \ast \mu_{\alpha} = \mu_{\alpha} \circ \pi_{2\alpha}.
$$
In the terms of transforms we arrive at:
$$
\Phi_{\pi_{2\alpha}}(t)^n = \Phi_{\lambda_{0,2n,\alpha}(\delta_1)}(t).
$$
It means that we are able to cumulate the Pareto distributions in the Kendall convolution algebra.

Now we construct random walk under the Kendall convolution with the unit step $\triangle X_{1}$ with distribution $\delta_1$.

\begin{thm}\label{thm:3}
The Markov process $\{X_n: n\in\mathbb{N}_0\}$, with $X_0=0$ and $ X_{1} \sim \delta_1$, based on the Kendall convolution $\triangle_{\alpha}$, has the following properties:
$$
X_2=\theta_1,\; X_{n+1} = X_n \cdot \theta_{n}^{Q_n} \quad a.e.,
$$
where $\theta_n$ is independent of $X_n$  for $n \geqslant 1$ and
$$P(Q_n=k|X_n) = \left\{ \begin{array}{lcl}
      \frac{1}{X_n^{\alpha}}, & \hbox{ for } & k=1; \\
      1-\frac{1}{X_n^{\alpha}}, & \hbox{ for } & k=0.
\end{array} \right.
$$
Moreover
\begin{eqnarray*}
  \lefteqn{
P_{n-1,n}(x,A):=P(X_n\in A|X_{n-1}=x)}\\
& & =\frac{1}{x^{\alpha}}P(x \theta_{n-2}\in
A)+\left(1-\frac{1}{x^{\alpha}}\right)\mathbf{1}_{A}(x)
\end{eqnarray*}
for every Borel set $A\subseteq[0,\infty)$.
\end{thm}
\noindent {\bf Proof.}
By the construction the transition probabilities are given by
\begin{eqnarray*}
\lefteqn{\hspace{-5mm}P(X_n\in A|X_{n-1}=x) = \delta_x \diamond \delta_1^{\diamond n-k}(A)} \\
 & = & \frac{1}{x^{\alpha}}P\Bigl(x \theta_{n-1}\in A\Bigr) +\Bigl(1-\frac{1}{x^{\alpha}} \Bigr)\mathbf{1}_A(x).
\end{eqnarray*}
We see that
\begin{eqnarray*}
  \lefteqn{\left(\lambda_{0,n,\alpha}(\delta_1)\right)(A)=
  \int_{1}^{\infty} \left(\delta_x \triangle_{\alpha} \delta_1\right)(A)\left(\lambda_{0,n-1,\alpha}(\delta_1)\right)(dx)=}\\
& & =\int_{1}^{\infty}\left(T_x\Bigl(\frac{1}{x^{\alpha} }\,\delta_1\triangle_{\alpha}\delta_1 + \Bigl(1- \frac{1}{x^{\alpha}}\Bigr)\delta_1\Bigr)\right)(A)\left(\lambda_{0,n-1, \alpha}(\delta_1)\right)(dx)\\
  & & =\int_{1}^{\infty}\left(\frac{1}{x^{\alpha}}P(x \theta_{n-1}\in A)+\left(1-\frac{1}{x^{\alpha}}\right)
  \mathbf{1}_A(x)\right)\left(\lambda_{0,n-1,\alpha}(\delta_1)\right)(dx).
\end{eqnarray*}
for all $A\in\mathcal{B}((0,\infty))$ and $n\geqslant 2$. By Corollary 5.11 in \cite{Kall} we can find such sequences $\{\theta_n\}$ and $\{Q_n\}$ that
$$
X_{n+1}=X_n \cdot \theta_{n}^{Q_n} \quad a.e.
$$
Finally, $X_n$ has distribution $\lambda_{0,n,\alpha}(\delta_1)$.
\qed

\vspace{2mm}

In the next theorem we give generalization of Theorem \ref{thm:3} for $\Delta X_1$ with any distribution $\nu$ concentrated on the positive half line.
\begin{thm}\label{thm:4}
The Markov process $\{X_n: n\in\mathbb{N}\}$ with $\Delta X_1 \sim \nu\in\mathcal{P_+}$ based on the Kendall convolution $\triangle_{\alpha}, \,\alpha>0$, has the following properties:
$$
X_0=0, \; X_1= \Delta X_1, \; X_{n+1}=\left(X_n \vee \triangle X_{n+1} \right)\cdot \theta_{n}^{Q_n} \quad a.e.,
$$
where $\theta_{n}$ is independent of $\left(X_n \vee \triangle X_{n+1} \right)$ for $n \geqslant 1$ and
$$
P(Q_n=k|X_n,\triangle X_{n+1}) = \left\{ \begin{array}{lcl}
      z(X_{n},\triangle X_{n+1})^{\alpha}, & \hbox{ for } & k=1; \\
      1-z(X_{n},\triangle X_{n+1})^{\alpha}, & \hbox{ for } & k=0.
\end{array} \right.
$$
for $z(x,y)=\frac{x\wedge y}{x\vee y}$.
\\
Moreover
\begin{eqnarray*}
  \lefteqn{
P_{n-1,n}(x,y,A):=P(X_n\in A|X_{n-1}=x,\triangle X_n=y)}\\
& & =z(x,y)^{\alpha}P(v(x,y) \cdot \theta_{n-2}\in
A)+\left(1-z(x,y)^{\alpha}\right)\pmb{1}_{A}(v(x,y))
\end{eqnarray*}
for $v(x,y)=x\vee y$ and every Borel set $A\subseteq[0,\infty)$.
\end{thm}
\noindent {\bf Proof.}
Similarly, as in the proof of Theorem \ref{thm:3} we can show that
\begin{eqnarray*}
\lefteqn{
  \left(\lambda_{0,n,\alpha}(\nu)\right)(A)=\left(\left(\lambda_{0,n-1, \alpha}(\nu)\right)\triangle_{\alpha}\nu\right)(A)}\\
  & & =
  \int\limits_0^{\infty} \int\limits_0^{\infty} \left(\delta_x \triangle_{\alpha} \delta_y\right)(A)\left(\lambda_{0,n-1, \alpha}(\nu)\right)(dx)\nu(dy).
\end{eqnarray*}
  Then substituting $v(x,y)=v$ and $z(x,y)=z$ we arrive at
\begin{eqnarray*}
\lefteqn{\int\limits_0^{\infty} \int\limits_0^{\infty}\left(T_v \left(z^{\alpha}\,(\delta_1\triangle_{\alpha}\delta_{z}) + \left(1- z^{\alpha}\right)\delta_1\right)\right)(A)\left(\lambda_{0,n-1, \alpha}(\nu)\right)(dx)\nu(dy)}\\
  & & =\int\limits_0^{\infty} \int\limits_0^{\infty}\left(z^{\alpha}P(v \theta_{n-1}\in A)+\left(1-z^{\alpha}\right)
  \mathbf{1}_A(v(x,y))\right)\left(\lambda_{0,n-1, \alpha}(\nu)\right)(dx)\nu(dy)
\end{eqnarray*}
for all $A\in\mathcal{B}((0,\infty))$ and $n\geqslant 2$. It yields by Corollary 5.11 in \cite{Kall} that
$$
X_{n+1}=v \left(X_n,\triangle X_{n+1}\right)\cdot \theta_{n}^{Q_n} \quad a.e.
$$
for $Q_n$ with desired distribution, the proper sequence $\{\theta_n\}$ and $X_n$ has distribution $\lambda_{0,n,\alpha}(\nu)$.
\qed

\vspace{2mm}

\par
One would think that the Markov chain $\{X_n: n\in\mathbb{N}_0\}$ with respect to a generalized convolution $\diamond$, with unit step distribution $\nu \in \mathcal{P}_+$, is a kind of the L\'evy process in the classical sense. The following two propositions give a negative answer to such a hypothesis.

\begin{prop}\label{prop:2}
Let $\triangle_1$ be the Kendall convolution and $\{X_n: n\in\mathbb{N}\}$ be a random walk under $\triangle_1$ with unit step distribution $\delta_1$. Then
$$
P(X_{k+1}-X_k<w)=1-\frac{2}{k+1}E\left((1+wY)^{-2}\right)
$$
for every $k\in\mathbb{N}$, where $Y$ has distribution $\beta_{3,k-1}$, which means that the increments of this chain are not stationary in time.
\end{prop}
\noindent {\bf Proof.}
Notice that by constructing of random walk given in Theorem \ref{thm:4}, we have
$$
P(X_{k+1}=X_k)=\int\limits_1^{\infty}P(Q_k=0|X_k=s)\left(\lambda_{0,k, \alpha}(\delta_1)\right)(ds)=\frac{k-1}{k+1}.
$$
The continuous part of the distribution of $\bigl(X_k ,X_{k+1}\bigr)$ has the weight $2/{(k+1)}$ and the density
$$
f(u,v)=\frac{(k+1)k(k-1)}{u^2 v^3}\; \Bigl(1-\frac{1}{u} \Bigr)^{k-2} \! \pmb{1}_{\{1\leqslant u\leqslant v\}}.
$$
Since $X_1\sim\delta_1$ and
\begin{eqnarray*}
P\Bigl\{X_{k+1}-X_k<w \Bigr\} & = & \frac{k-1}{k+1} +  \int\limits_1^{\infty} \int\limits_{u}^{u+w}  f(u,v) dv du \\
& = & 1 - k(k-1) \,\int\limits_1^{\infty}\frac{({u}-1)^{k-2}}{u^k\,(u+w)^2}\, du,
\end{eqnarray*}
thus the Markov process $\{X_n: n\in\mathbb{N}_0\}$ does not have stationary increments.
\qed
\begin{prop}\label{prop:3}
The increments of the Markov chain $\{X_n: n\in\mathbb{N}\}$ with respect to the Kendall convolution $\triangle_1$, with unit step $\delta_1$, are not independent. In particular, for each $k\in\mathbb{N}$ the random variables $X_k$ and $X_{k+1}-X_k$ are not independent.
\end{prop}
\noindent {\bf Proof.}
By simple computation we arrive at
\begin{eqnarray*}
\lefteqn{
P(X_{k+1}-X_k<w,X_k<z)=P(X_{k+1}-X_k<w,X_k<z|X_{k+1}=X_k)}\\
& & \cdot P(X_{k+1}=X_k)+P(X_{k+1}-X_k<w,X_k<z,X_{k+1}>X_k)\\
& & = \frac{k-1}{k+1}\left(\lambda_{0,k,\alpha}(\delta_1)\right)([0,z))+ \int\limits_1^{z}\int\limits_{u}^{u+w} f(u,v) dv du \\
& & = \frac{k-1}{k+1}\left(\lambda_{0,k,\alpha}(\delta_1)\right)([0,z))+ \frac{2}{k+1}\left[1-B(3,k-1,z^{-1})\right]\\
& & -\frac{2}{k+1}\int\limits_{1/z}^{1}\frac{k(k-1)}{(1+wy)^2}(1-y)^{k-2} dy  = \frac{k-1}{k+1}\left(\lambda_{0,k,\alpha}(\delta_1)\right)([0,z))\\
& & + \frac{2}{k+1}\left[1-B(3,k-1,z^{-1})\right]-\frac{2}{k+1}E\left(\frac{1}{(1 +wY)^2}\pmb{1}_{\{\omega: Y(\omega)>1/w\}}\right),
\end{eqnarray*}
where
$Y\sim \beta_{3,k-1}$ and $B(3,k-1,z^{-1})$ is the function given in Corollary \ref{cor:2}. \\
In particular
\begin{eqnarray*}
\lefteqn{
P(X_3-X_2<w, X_2<z)= 1-\frac{1}{3z^2}\left(1+\frac{2}{z}\right)}\\
& & -\frac{2}{w^3}\left[ln\left(\frac{w+z}{z(1+w)}\right)^2+\frac{w(z-1)}{z (w+z)(1+w)}\left(z-(w+z)(1+w)\right)\right]
\end{eqnarray*}
for $w>0,z>1$ and $P(X_2<z)=\left(1-\frac{1}{z^2}\right)$ for $z>1$. Moreover
$$
P(X_3-X_2<w)=\frac{2}{3}-2\left[\frac{1}{w^2}-\frac{w}{1+w}+\frac{3}{w^2 (1+w)^2}-\frac{1}{w^3}ln\left(1+w\right)^2\right]
$$
for $w>0$, which yields that $X_2$ and $X_3-X_2$ are not independent.
\qed
\vspace{2ex}

Notice that in the same manner we can construct a family of random walks $\{Z_n^{(k)}:n\in\mathbb{N}\}$ under the Kendall convolution such that
$$
Z_n^{(k)}:=X_{kn}
$$
for every $k\in\mathbb{N}$. For every fixed $k$ the unit step of this random walk has distribution $\lambda_{0,k,\alpha}(\nu)$ given in Proposition \ref{prop:1}. In particular, one can prove that probability measures which belong to the family $\{\lambda_{0,k,\alpha}(\delta_1):k\in\mathbb{N}\}$ are heavy tailed and $\lambda_{0,2,\alpha}(\delta_1)=\pi_{2\alpha}$.

\medskip
In a similar way we construct random walk under the Kendall weak generalized convolution $\otimes_{\mu_{\alpha}}$, where $\alpha\in(0,1]$. Let $Y,(Y_i)_{i\in\mathbb{N}_0}$ be a sequence of $i.i.d.$ random variables with weakly stable distribution $\mu_{\alpha}$,  $(\Delta \widetilde{X}_i)_{i\in\mathbb{N}_0}$ be a sequence of $i.i.d$ random variables with distribution $\nu$, which are concentrated on the real line. Let
$(\widetilde{\theta}_i)_{i\in\mathbb{N}_0}$ be a sequence of $i.i.d$ random variables with distribution $\widetilde{\pi}_{2\alpha}$ such that $\widetilde{\theta}_n$ is independent of $\widetilde{X}_n$ and $\Delta \widetilde{X}_{n+1}$ for all $n\in\mathbb{N}$. Moreover, let the sequences $Y,(Y_i)_{i\in\mathbb{N}_0}$ and $(\widetilde{\theta}_i)_{i\in\mathbb{N}_0}$ be independent. Then $\widetilde{X}_n$ denotes the position of moving particle at the n-th step such that the unit step has distribution $\nu$ and $\widetilde{X}_n$ has distribution being symmetrization of $\lambda_{0,n,\alpha}(\nu)$, i.e.
$$
\widetilde{\lambda}_{0,n,\alpha}(\nu):=\nu^{\otimes_{\mu_{\alpha}}n},
$$
where $\widetilde{\lambda}_{0,n,\alpha}(\nu)$ is the probability measure on the real line and $\widetilde{\lambda}_{0,0,\alpha}(\nu)=\delta_0$.
Just as in the case of random walk under the Kendall generalized convolution we can get the series of dual results for random walk with respect to $\otimes_{\mu_{\alpha}}$:
\begin{lem}\label{lem:1}
For each natural number $n\geqslant 2$ and $\alpha\in(0,1]$ probability measure $\widetilde{\lambda}_{0,n,\alpha}(\widetilde{\delta}_1)$ has the density
$$
\left(\widetilde{\lambda}_{0,n,\alpha}(\widetilde{\delta}_1)\right)(dx) =\frac{\alpha n(n-1)}{2|x|^{2\alpha+1}} \left(1-\frac{1}{|x|^{\alpha}}\right)^{n-2} \pmb{1}_{[1,\infty)}(|x|)dx.
$$
\end{lem}
The above Lemma is a modification of Corollary \ref{cor:4} except that here we take the characteristic function of $\mu_{\alpha}$ as the kernel of the corresponding homomorphism. In the next theorem we construct the Markov process with distribution $\widetilde{\lambda}_{0,n,\alpha}(\nu)$ at the n-th step.

\begin{thm}\label{thm:5}
The Markov process $\{\widetilde{X}_n: n\in\mathbb{N}_0\}$ with $\triangle \widetilde{X}_1 \sim \nu$ has the following properties:
 $$
\widetilde{X}_0=0, \; \widetilde{X}_1= \Delta \widetilde{X}_1, \; \widetilde{X}_{n+1}= v\left(|\widetilde{X}_n|, |\triangle \widetilde{X}_{n+1}| \right)\cdot u\left(\widetilde{X}_n, \triangle \widetilde{X}_{n+1} \right) \cdot \widetilde{\theta}_{n}^{\widetilde{Q}_n} \quad a.e.
$$
for $n \geqslant 1$, where $\theta_{n}$ is independent of $v\left(|\widetilde{X}_n|, |\triangle \widetilde{X}_{n+1}| \right)\cdot u\left(\widetilde{X}_n, \triangle \widetilde{X}_{n+1} \right)$,
$v(x,y)=x\vee y, \; z(x,y)=\frac{x\wedge y}{x\vee y}$,
$$
u(x,y) = \left\{ \begin{array}{lcl}
      sgn(x), & \hbox{ for } & |x| \geqslant |y|, \\
      sgn(y), & \hbox{ for } & |y| \geqslant |x|,
\end{array} \right.
$$
and
$$
P(\widetilde{Q}_n=k|\widetilde{X}_n,\triangle \widetilde{X}_{n+1}) = \left\{ \begin{array}{lcl}
      \left(z(|\widetilde{X}_n|,|\triangle \widetilde{X}_{n+1}|)\right)^{\alpha}, & \hbox{ for } & k=1; \\
      1-\left(z(|\widetilde{X}_n|,|\triangle \widetilde{X}_{n+1}|)\right)^{\alpha}, & \hbox{ for } & k=0.
\end{array} \right.
$$
Moreover $\widetilde{X}_n$ has distribution $\widetilde{\lambda}_{0,n,\alpha}(\nu)$ for every $n\in\mathbb{N}$ and
\begin{eqnarray*}
  \lefteqn{
P_{n-1,n}(x,y,A):=P(\widetilde{X}_n\in A|\widetilde{X}_{n-1}=x,\triangle \widetilde{X}_n=y)}\\
& & = \left(z(|x|,|y|)\right)^{\alpha} P( u(x,y)\cdot v(|x|,|y|) \cdot \widetilde{\theta}_{n-2}\in
A)\\
& & + \frac{1}{2}\left(1-\left(z(|x|,|y|)\right)^{\alpha}\right)\pmb{1}_{\widetilde{A}}(u(x,y)\cdot v(|x|,|y|))
\end{eqnarray*}
for every Borel set $A\subseteq \mathbb{R}$, where $\widetilde{A}$ is symmetrization of $A$.
\end{thm}
\noindent {\bf Proof.}
In order to prove this theorem it is sufficient to follow the proof of Theorem \ref{thm:4} substituting $\otimes_{\mu_{\alpha}}$ instead of $\triangle_{\alpha}$ and integrating over $\mathbb{R}$ but using the probability kernel of weak Kendall convolution. It is easy to see that for every Borel set $A\in\mathcal{B}(\mathbb{R})$ we have
\begin{eqnarray*}
\lefteqn{
  \left(\widetilde{\lambda}_{0,n,\alpha}(\nu)\right)(A)= \int\limits_{\mathbb{R}} \int\limits_{\mathbb{R}}\Big(\left(z(|x|,|y|)\right)^{\alpha} P( u(x,y)\cdot v(|x|,|y|) \cdot \widetilde{\theta}_{n-2}\in A)}\\
& & + \left(1-\left(z(|x|,|y|)\right)^{\alpha}\right)\widetilde{\delta}_A( u(x,y) \cdot v(|x|,|y|)) \Big)\nu(dy)
\left(\widetilde{\lambda}_{0,n-1, \alpha}(\nu)\right)(dx).
\end{eqnarray*}
\qed

In particular for $\nu=\widetilde{\delta}_1$ we have the following construction.
\begin{cor}\label{cor:5}
The Markov process $\{\widetilde{X}_n: n\in\mathbb{N}\}$ with $\triangle \widetilde{X}_1 \sim \widetilde{\delta}_1$ has the following properties:
 $$
\widetilde{X}_0=0, \; \widetilde{X}_1= \Delta \widetilde{X}_1, \; \widetilde{X}_2=\widetilde{\theta}_1,\; \widetilde{X}_{n+1}=\widetilde{X}_n \cdot \widetilde{\theta}_{n}^{\widetilde{Q}_n} \quad a.e.
$$
for $n\geqslant 2$, where $\theta_{n}$ is independent of $\widetilde{X}_n$ and
$$P(\widetilde{Q}_n=k|\widetilde{X}_n) = \left\{ \begin{array}{lcl}
      \frac{1}{|\widetilde{X}_n|^{\alpha}}, & \hbox{ for } & k=1; \\
      1-\frac{1}{|\widetilde{X}_n|^{\alpha}}, & \hbox{ for } & k=0.
\end{array} \right.
$$
Moreover $\widetilde{X}_n$ has distribution $\widetilde{\lambda}_{0,n,\alpha}(\widetilde{\delta}_1)$ for every $n\in\mathbb{N}$ and
\begin{eqnarray*}
  \lefteqn{
\widetilde{P}_{n-1,n}(x,A):=P(\widetilde{X}_n\in A|\widetilde{X}_{n-1}=x)}\\
& & =\frac{1}{|x|^{\alpha}}P(x \widetilde{\theta}_{n-2}\in
A) + \frac{1}{2} \left(1-\frac{1}{|x|^{\alpha}}\right)\mathbf{1}_{\widetilde{A}}(x)
\end{eqnarray*}
for every Borel set $A\in\mathcal{B}(\mathbb{R})$, where $\widetilde{A}$ is symmetrization of $A$.
\end{cor}
Notice that using the sequence $(Y_i)_{i\in\mathbb{N}}$ we have constructed also a random walk in the usual sense $\{\widetilde{S}_n:n\in\mathbb{N}_0\}$ associated with the random walk $\{\widetilde{X}_n:n\in\mathbb{N}_0\}$ under weak Kendall convolution $\otimes{\mu_{\alpha}}$ with unit step $\triangle \widetilde{X}_1 \sim \widetilde{\delta}_1$. The unit step $\widetilde{S}_1$ has distribution $\mu_{\alpha}$ and
$$
\widetilde{S}_0=0, \; \widetilde{S}_n=Y_1+Y_2+\cdots+Y_n \quad a.e.
$$
The relation between $\{\widetilde{S}_n:n\in\mathbb{N}_0\}$  and $\{\widetilde{X}_n:n\in\mathbb{N}_0\}$ is given by the following distribution equation
$$
\widetilde{S}_n \stackrel{d}{=} Y \widetilde{X}_n, \quad \hbox{where}  \; Y \sim \mu_{\alpha},
$$
for every $n\in\mathbb{N}_0$.
\\
In the next few lemmas we give some properties for the random walk $\{\widetilde{S}_n:n\in\mathbb{N}_0\}$ associated with $\{\widetilde{X}_n: n\in\mathbb{N}_0\}$, where $\triangle\widetilde{X}_1\sim\widetilde{\delta}_1$. In particular, we have the following relation between $\mu_{\alpha}$ and $\mu_1$.
\begin{lem}\label{lem:2}
Let $0<\alpha\leqslant 1$ and $\mu_{\alpha}$ be weakly stable probability measure, which induces the Kendall convolution. Then
$$
\mu_{\alpha} = \mu_1 \circ \left( \alpha \widetilde{\delta}_1 + (1-\alpha) \widetilde{\pi}_{\alpha}\right)
$$
\end{lem}
\noindent{\bf Proof.} In order to find a probability measure $\nu \in \mathcal{P}$ such that $\mu_{\alpha}= \mu_1 \circ \nu$ we have to solve the integral equation:
$$
\left(1-|t|^{\alpha}\right)_+ = \int\limits_{\mathbb{R}} \left(1-|ts|\right)_+ \nu(ds).
$$
This equation can be solved in the same manner as in the proof of Proposition \ref{prop:1}, however, it would be much simpler to check that for the measure $\nu = \alpha \widetilde{\delta}_1 + (1-\alpha) \widetilde{\pi}_{\alpha}$ the desired equality holds.
\qed
\begin{prop}\label{prop:2}
Let $0<\alpha\leqslant 1$ and $(Y_i)_{i \in\mathbb{N}_0}$ be the sequence of independent identically distributed random variables with distribution $\mu_{\alpha}$. Then $\widetilde{S}_n=Y_1+Y_2+\cdots+Y_n$ has distribution
$$
\mu_{\alpha}^{\ast n} = \left(\mu_1 \circ \left( \alpha \widetilde{\delta}_1 + (1-\alpha) \widetilde{\pi}_{\alpha}\right)\right)^{\ast n} = \mu_1 \circ \left( \alpha \widetilde{\delta}_1 + (1-\alpha) \widetilde{\pi}_{\alpha}\right)^{\otimes_{\mu_1} n},
$$
where
\begin{eqnarray*}
\lefteqn{
\left( \alpha \widetilde{\delta}_1 + (1-\alpha) \widetilde{\pi}_{\alpha}\right)^{\otimes_{\mu_1} n}(dx)}
\\
& = & \frac{\alpha n}{2|x|^{\alpha+1}}\left(1-\frac{1}{|x|^{\alpha}}\right)^{n-2} \left[1-\alpha + \frac{\alpha n-1}{|x|^{\alpha}}\right] \pmb{1}_{(1,\infty)}(|x|) \, dx
\end{eqnarray*}
\end{prop}
\noindent{\bf Proof.}
By Lemma \ref{lem:2} the characteristic function of measure $\mu_1 \circ \left( \alpha \widetilde{\delta}_1 + (1-\alpha) \widetilde{\pi}_{\alpha}\right)$ is given by the formula:
$$
G(1/t):=\alpha\widehat{\mu}_1(t) + (1-\alpha)\widehat{\mu_1\circ \pi_{\alpha}}(t) = \widehat{\mu}_{\alpha}(t)
$$
Since
$$
\left(G(1/t)\right)^n=\int_0^{\infty} (1-tx)_+ F_n(dx)=t\int\limits_0^{1/t} F_n(x) dx,
$$
where $F_n$ denotes the distribution function of $\left( \alpha \widetilde{\delta}_1 + (1-\alpha) \widetilde{\pi}_{\alpha}\right)^{\otimes_{\mu_1} n}$, and substituting $x:=1/t$,  we arrive at
$$
F_n(x)=\frac{d}{dx}\left[xG(x)^n\right],
$$
which leads to the explicit formula for $F_n$. \qed

\vspace{2mm}

The above proposition says that, by weak stability of $\mu_{\alpha}, \,\alpha \in (0,1]$, we can consider random walk under the Kendall convolution $\otimes_{\mu_1}$ with unit step $\triangle X_1 \sim \alpha \widetilde{\delta}_1 + (1-\alpha) \widetilde{\pi}_{\alpha}$ instead of random walk under $\otimes_{\mu_{\alpha}}$ with unit step with distribution $\widetilde{\delta}_1$. By this property we see that random walk in the usual sense with unit step $\mu_{\alpha}$ is also associated with the random walk under the Kendall convolution $\otimes_{\mu_1}$.
Similar property we have also for random walk with unit step with distribution $\nu\in\mathcal{P}$.
\begin{rem}\label{rem:1}
Let $\nu \in \mathcal{P}$ and $\alpha\in(0,1]$. By Lemma \ref{lem:2} and weak stability of measures $\mu_{\alpha}$ we have
\begin{eqnarray*}
\lefteqn{
\mu_1 \circ \left( \alpha \widetilde{\delta}_1 + (1-\alpha) \widetilde{\pi}_{\alpha}\right)\circ \nu^{\otimes_{\mu_{\alpha}} n} } \\
& & = \left(\mu_1 \circ \left( \alpha \widetilde{\delta}_1 + (1-\alpha) \widetilde{\pi}_{\alpha}\right) \circ \nu\right)^{\ast n}= \mu_1 \circ \left( \left( \alpha \widetilde{\delta}_1 + (1-\alpha) \widetilde{\pi}_{\alpha}\right)  \circ \nu \right)^{\otimes_{\mu_{1}} n}\!\!,
\end{eqnarray*}
which implies that
$$
\left( \alpha \widetilde{\delta}_1 + (1-\alpha) \widetilde{\pi}_{\alpha}\right)\circ \nu^{\otimes_{\mu_{\alpha}} n} = \left(\left( \alpha \widetilde{\delta}_1 + (1-\alpha) \widetilde{\pi}_{\alpha}\right) \circ \nu \right)^{\otimes_{\mu_{1}} n},
$$
\end{rem}
It is worth noticing that for the random walk with unit step $\mu_1$ we have the following recurrence relation.
\begin{lem}\label{lem:3}
For every natural number $n\geqslant 3$ we have
\begin{eqnarray*}
\lefteqn{
\mu_{1}(dx)  = \frac{1}{\pi x^2} \left(1-cosx\right) dx} \\
\mu_{1}^{\ast n}(dx) & = & \frac{n}{\pi x^2} dx - \frac{(n-1)n}{\pi x^2} \mu_{1}^{\ast n-2}(dx)
\end{eqnarray*}
\end{lem}
\noindent {\bf Proof.} Denoting by $g_n$ the density of the measure $\mu_1^{\ast n}$ and applying the Fourier inverse transform we obtain
$$
g_1(x)=\frac{1}{2\pi}\int\limits_{0}^{1} cos(tx)(1-t)dt,
$$
and
$$
g_n(x)=\frac{1}{\pi}\int\limits_{0}^{1} cos(tx)(1-t)^n dt.
$$
Now it is sufficient to use integration by parts. \qed

\vspace{2mm}

\begin{lem}\label{lem:4}
Let  $0<\alpha\leqslant 1$  and  $\otimes_{\mu_{\alpha}}$  be the weak Kendall convolution. Then random walk $\{\widetilde{X}_n: n\in\mathbb{N}\}$ with unit step $\triangle\widetilde{X}_1 \sim \widetilde{\delta}_1$ is not recurrent.
\end{lem}
\noindent {\bf Proof.} By Corollary \ref{cor:4} we have
$$
\lambda_{0,n,\alpha}(\delta_1)((0,x])=\left(1+\frac{n-1}{x^{\alpha}} \right)\left(1-\frac{1}{x^{\alpha}}\right)_+^{n-1}
$$
Since
$$
\left|\widetilde{\lambda}_{0,n,\alpha}(\widetilde{\delta}_1) \right|((-\infty,x])=\lambda_{0,n,\alpha}(\delta_1)((0,x])
$$
and $|\widetilde{X}_1| \sim \delta_1$ we arrive at
$$
\sum\limits_{n=1}^{\infty}P(|\widetilde{X}_n|<x)=x^{\alpha} (2-x^{-\alpha}) \mathbf{1}_{[1,\infty)}(x) <\infty.
$$
By the Borel-Cantelli Lemma we obtain
$$
P\left(\limsup_{n\to\infty}\left\{|\widetilde{X}_n|<x\right\}\right)=0,
$$
which implies that $\{\widetilde{X}_n:\; n\in\mathbb{N}_0\}$ is not recurrent. \qed
\\
Now we present the result describing the magnitude of the fluctuations for random walk under weak Kendall convolution.

\begin{prop}\label{prop:3}
For every $r>\frac{1}{2}$ and random walk under the Kendall convolution $\{\widetilde{X}_n:n\in\mathbb{N}_0\}$ with unit step $\triangle \widetilde{X}_1 \sim \widetilde{\delta}_1$ we have
$$
\mathbf{P} \Biggl( \bigcup_{n=1}^{\infty} \bigcap_{k=n}^{\infty} \left\{ \bigl|\widetilde{X}_n \bigr|^{\alpha} \leqslant \frac{n^{r+1}}{\ln(n)} \right\} \Biggr) = 1.
$$
\end{prop}
\noindent {\bf Proof.}
To see this it is sufficient to notice that
$$
\mathbf{P} \left( |\widetilde{X}_n|^{\alpha} > \frac{n^{r+1}}{\ln(n)} \right) = 1 - \left(1+\frac{n-1}{n^{r+1}}\ln n \right)\left(1- n^{-r-1} \ln n \right)_+^{n-1}.
$$
Let $A_n = \left\{ |\widetilde{X}_n|^{\alpha} > \frac{n^{r+1}}{\ln(n)}\right\}$.
It is a matter of laborious but straightforward calculations to show that
$$
\lim_{n\rightarrow\infty} \frac{\mathbf{P}(A_n)}{n^{-2r} (\ln n)^2} = 1.
$$
Moreover for $2r-1>0$ we have
$$
\int_1^{\infty} x^{-2r} (\ln x)^2 dx = \int_0^{\infty} u^2\, e^{(1-2r)u}\, du = \frac{2}{(2r-1)^{3}} < \infty
$$
and we arrive at the assertion.
\qed

\begin{rem}\label{rem:2}
 Since $\mathbb{E} (|\widetilde{X}_n|^{\alpha})=n$ by the Tchebyshev inequality we also have that
$$
\mathbf{P} \left\{ |\widetilde{X}_n|^{\alpha} \geqslant \frac{n^{2}}{\ln(n)} \right\} \leqslant \frac{\mathbb{E} (|\widetilde{X}_n|^{\alpha})}{n^2} \ln n = \frac{\ln n}{n} \rightarrow 0.
$$
for $n \to \infty$.
\end{rem}

\end{document}